\documentclass[10pt,twoside]{article}
\usepackage{amsmath,amssymb}
\usepackage{Latex-document}

\newtheorem{lemma}{Lemma}

\newtheorem{theorem}[lemma]{Theorem}
\newtheorem{corollary}[lemma]{Corollary}
\newtheorem{exercise}[lemma]{Exercise}
\newtheorem{example}[lemma]{Example}
\newtheorem{definition}[lemma]{Definition}

\def\Bscr{{\cal B}}
\def\Cscr{{\cal C}}
\def\Dscr{{\cal D}}

\def\Fscr{{\cal F}}

\def\Iscr{{\cal I}}
\def\Jscr{{\cal J}}

\def\Lscr{{\cal L}}

\def\Oscr{{\cal O}}
\def\Pscr{{\cal P}}

\def\Rscr{{\cal R}}
\def\Sscr{{\cal S}}

\def\collon{{\hskip -1.5pt}:{\hskip -1.5pt}}
\def\Skl{\operatorname{Skl}}
\def\HH{\operatorname{H}}
\def\Num{\operatorname{Num}}
\def\GKdim{\operatorname{GKdim}}
\def\pone{\mathbb P^1}
\def\too{\longrightarrow}

\newcommand{\Proj}{\operatorname{Proj}}

\def\mod{\operatorname{mod}}

\def\qgr{\operatorname{qgr}}
\def\gr{\operatorname{gr}}
\def\coh{\operatorname{coh}}

\def\Ext{\operatorname {Ext}}
\def\Hom{\operatorname {Hom}}
\def\Id{\operatorname {Id}}
\def\gldim{\operatorname {gl\,dim}}

\DeclareMathOperator{\tors}{tors}
\DeclareMathOperator{\Aut}{Aut}

\title{\bf Noncommutative Projective Geometry\thanks{The author is supported in part by the National Science Foundation
under grant DMS-9801148.}\vskip 6mm}

\author{J. T. Stafford\vspace*{-0.5cm}\thanks{Department of Mathematics,
The University of Michigan, Ann Arbor, MI 48109-1109, USA. E-mail: jts@umich.edu}}

\date{\vspace{-8mm}}

\begin{document}

\maketitle

\thispagestyle{first} \setcounter{page}{93}

\begin{abstract}

\vskip 3mm

This  article describes recent applications of algebraic geometry to  noncommutative algebra. These techniques
have been particularly successful in describing graded algebras of small dimension.

\vskip 4.5mm


\end{abstract}

\vskip 12mm

\section{Introduction} \label{sec0}\setzero
\vskip-5mm \hspace{5mm}

In recent years a surprising number of significant insights and results in
noncommutative algebra have  been obtained by using the global techniques
 of projective algebraic geometry. This article will survey some
of these results.

The classical approach to projective geometry, where one  relates a
commutative graded domain $C$ to the associated variety $X=\Proj C$ of
homogeneous, nonirrelevant prime ideals, does not generalize well to the
noncommutative  situation, simply because noncommutative algebras do not have
enough ideals. However, there is a second approach, based on a classic theorem
of Serre: If $C$ is generated in degree one, then the categories
$\coh(X)$  of coherent sheaves on $X$ and  $\qgr C$ of finitely generated
graded $C$-modules  modulo torsion are equivalent.

Surprisingly, noncommutative analogues of this idea work very well  and have
lead to a  number of deep results.  There are two strands to
this approach. First, since  $X$ can be reconstructed from $\coh(X)$
\cite{rosenberg1} we will regard $\coh(X)$ rather than  $X$ as the variety
since this is what generalizes. Thus, given a  noncommutative graded
$k$-algebra $R=\bigoplus R_i$ generated in degree one  we will consider $\qgr
R$ as the corresponding ``noncommutative variety'' (the formal definitions will
be given in a moment). In particular, we will regard $\qgr R$ as a
noncommutative curve, respectively surface, if $\dim_kR_i$ grows linearly,
respectively quadratically. This analogy works  well,  since there are many
situations in which one can pass back and forth  between $R$ and $\qgr R$
\cite{AZ} and, moreover, substantial geometric  techniques can be applied to
study $\qgr R$. A survey of this approach may be found in \cite{SV}.

The second strand is more concrete. In order to use  algebraic geometry to
study noncommutative algebras we need to be able to create honest varieties
from those algebras. This is frequently possible  and such an approach will
form the basis of this survey. Once again, the idea is simple: when $R$ is
commutative, the points  of $\Proj R$  correspond to the graded factor modules
$M = R/I = \bigoplus _{i\geq 0} M_i$  for which $\dim_k M_i=1$ for all $i$.
These modules are still defined when $R$  is noncommutative and are called
point modules. In many circumstances  the set of all such modules is
parametrized by a commutative scheme and that scheme controls the structure of
$R$.

This article surveys significant applications of this idea. Notably:

\begin{itemize}
\item If $R=\bigoplus R_i$ is a domain such that $\dim_k R_i$ grows linearly,
then $\qgr R\simeq \coh(X)$ for a curve $X$ and $R$ can be reconstructed
from data on $X$.  Thus, noncommutative curves are commutative
(see Section~4).

\item The noncommutative analogues $\qgr R$ of the projective plane can be
classified. In this case, the
point modules are parametrized by either $\mathbb P^2$ (in which case
$\qgr R\simeq \mathbb P^2$)
or by an cubic curve $E\subset \mathbb P^2$, in which case data on
$E$  determines $R$ (see Section~2).

\item For  strongly noetherian rings, as defined in
Section~5, the point modules are always
parametrized by a projective scheme. However there exist many
noetherian algebras $R$ for which  no such parametrization exists. This  has
interesting consequences for the classification of noncommutative surfaces.
\end{itemize}

We now make precise the definitions that will hold throughout this article.
All rings will be algebras over a fixed, algebraically closed base field $k$
(although most of the results actually hold for arbitrary fields). A
$k$-algebra $R$ is called {\em connected graded} ({\em cg}) if $R$ is a
finitely generated $\mathbb N$-graded $k$-algebra $R=\bigoplus_{i\geq 0} R_i $
with $R_0=k$. Note that this forces $\dim_k R_i<\infty$ for all $i$. Usually,
we will assume that $R$ is  generated in degree one in the sense that $R$ is
generated by $R_1$ as a $k$-algebra.  If $R=\bigoplus_{i\in \mathbb N}R_i$ is a
right noetherian cg ring then  define $\gr R$ to be the category of finitely
generated, $\mathbb Z$-graded right $R$-modules, with morphisms being graded
homomorphisms of degree zero. Define  the torsion subcategory,  $\tors R$, to
be the full subcategory of $\gr R$  generated by the  finite dimensional
modules and  write  $\qgr R=\gr R /\tors R$.  We write $\pi$ for the canonical
morphism $\gr R \to \qgr R$ and set $\Rscr = \pi(R)$.

 One can---and often should---work more generally with all graded $R$-modules
 and all quasi-coherent sheaves of $\Oscr_X$-modules,
  but two categories are enough.

In order to measure the growth of an algebra we use the following  dimension
function: For a cg ring $R=\bigoplus_{i\geq 0} R_i $, the  {\it
Gelfand-Kirillov dimension} of $R$ is defined to be  $ \GKdim R \;=\; \inf
\left\{\alpha\in \mathbb R : \dim_k (\sum_{i=0}^n R_i) \leq n^\alpha \text{ for
all } n\gg 0\right\}. $  Basic facts about this dimension can be found in
\cite{KL}.  If $R$ is a commutative cg algebra then  $\GKdim R$  equals the
Krull dimension of $R$ and hence equals $\dim \Proj R+1$. Thus a
noncommutative curve, respectively surface, will more formally
 be defined as  $\qgr R$  for
a cg algebra $R$ with $\GKdim R = 2$, respectively $3$.

\section{Historical background} \label{sec1.5} \setzero
\vskip-5mm \hspace{5mm}

We begin with a historical introduction to the subject. It really started with
the  work of Artin and Schelter \cite{AS} who attempted to classify the
noncommutative analogues $R$ of  a polynomial ring in three variables (and
therefore of $\mathbb P^2$).  The first problem is one of definition. A
``noncommutative polynomial ring'' should obviously be a cg ring of finite
global dimension, but this is too general, since it includes the free algebra.
One can circumvent this problem by requiring that $\dim_kR_i$ grows
polynomially, but  this still does not exclude unpleasant rings like
$k\{x,y\}/(xy)$ that has global dimension two but is neither noetherian nor a
domain.  The solution is to impose a Gorenstein condition and this leads to the
following definition:

\begin{definition}\label{ASreg}
 A cg algebra $R$ is called {\em AS-regular of dimension $d$}
 if  $\gldim R=d$,  $\GKdim R<\infty$
and  $R$ is AS-Gorenstein; that is,
$\Ext^i(k,R) =  0$ for $i\not=d$ but $\Ext^d(k,R)=k$,
 up to a shift of degree.
\end{definition}

One advantage with the Gorenstein hypothesis,  for AS-regular rings of dimension
$3$, is that the projective resolution of $k$ is forced to be of the form
$$0\too R\too R^{(n)} \too R^{(n)} \too R\too k\too 0$$ for some $n$
and, as Artin and Schelter show in \cite{AS},
 this  gives  strong information on the
Hilbert series and hence the defining relations of $R$.  In the process they
constructed one class of algebras that they were unable to analyse:

\begin{example} \label{skl} \rm
The three-dimensional Sklyanin algebra is the algebra
$$\Skl_3=\Skl_3(a,b,c)=
k\{x_0,x_1,x_2\}/(ax_ix_{i+1} + bx_{i+1}x_i + cx_{i+2}^2 :
i\in \mathbb Z_3),$$
where $(a,b,c)\in \mathbb P^2\smallsetminus F$, for a (known) set $F$.
\end{example}

The original Sklyanin algebra $\Skl_4$  is a $4$-dimensional analogue of
$\Skl_3$ discovered in  \cite{Skly1}. Independently of \cite{AS},
Odesskii and Feigin  \cite{OF2} constructed  analogues of $\Skl_4$ in all
dimensions and coined the name Sklyanin algebra.
 See \cite{connes} for applications of Sklyanin algebras to another
 version of noncommutative geometry.

In retrospect the reason  $\Skl_3$ is hard to analyse is because it depends
upon an elliptic curve and so a more geometric approach is required. This
approach came in \cite{ATV1} and depended upon the following simple idea.
Assume that $R$ is a cg algebra that is generated in degree one.  Define a {\it
point module} to be a cyclic  graded (right) $R$-module $M=\bigoplus_{i\geq 0}
M_i$ such that $\dim_kM_i=1$ for all $i\geq 0$.  The notation is justified by
the fact that, if $R$ were commutative, then  such a point module $M$ would be
isomorphic to $k[x]$ and hence equal to the homogeneous coordinate ring  of a
point in $\Proj R$. Point modules are easy to analyse geometrically  and  this
provides an avenue for using geometry in the study of cg rings.

We will illustrate this approach for  $S=\Skl_3$. Given a point module
$M=\bigoplus M_i$ write $M_i=m_ik$ for some $m_i\in M_i$ and suppose that the
module structure is defined by $m_ix_j = \lambda_{ij}m_{i+1}$ for some
$\lambda_{ij}\in k$.  If $f=\sum f_{ij}x_ix_j$ is one of the relations for $S$,
then necessarily $m_0f = (\sum f_{ij}\lambda_{0i}\lambda_{1j})m_2$, whence
$\sum f_{ij}\lambda_{0i}\lambda_{1j}=0$. This defines a subvariety
$\Gamma\subseteq \mathbb P(S_1^*)\times \mathbb P(S_1^*)=\mathbb P^2\times \mathbb
P^2$ and clearly $\Gamma$ parametrizes the {\it truncated point modules of
length three}: cyclic $R$-modules $M=M_0\oplus M_1\oplus M_2$ with
$\dim M_i =1$ for $0\leq i\leq 2$. A simple computation  (see
\cite[Section~3]{ATV1} or \cite[Section~8]{SV}) shows that $\Gamma$ is actually
the  graph of an automorphism $\sigma$ of an elliptic curve $E\subset \mathbb
P^2$. It follows easily  that $\Gamma$ also parametrizes the point
modules. As a morphism of point modules, $\sigma$ is nothing more than the
shift functor $M=\bigoplus M_i \mapsto M_{\geq 1}[1] = M_1\oplus
M_2\oplus\cdots.$

The next question is how to use $E$ and $\sigma$ to understand $\Skl_3$.
Fortunately, one can create a noncommutative algebra from this data that is
closely connected to  $\Skl_3$. This is the {\it twisted homogeneous coordinate
ring} of $E$  and is defined as follows. Let $X$ be a $k$-scheme, with a line
bundle $\Lscr$ and automorphism $\sigma$. Set $\Lscr_n =  \Lscr\otimes
\Lscr^\sigma\otimes\cdots\otimes \Lscr^{\sigma^{n-1}}$, where
$\Lscr^\tau=\tau^*\Lscr$ denotes the pull-back of $\Lscr$ along an automorphism
$\tau$. Then  the {\it twisted homogeneous coordinate ring} is defined to be
the graded vector space $B=B(X,\Lscr,\sigma)=k+\bigoplus_{n\geq 1} B_n$ where
$B_n=\HH^0(X,\Lscr_n)$. The multiplication on $B=B(Y,\Lscr,\sigma)$ is
defined by the   natural map
$$ \begin{array}{rl}
B_n\otimes_k B_m
&\cong \,\,
\HH^0({{X}}, \Lscr_n)\otimes_k \sigma^n\HH^0({{X}}, \Lscr_m)
\\ \noalign{\vskip 5pt}
&\cong\,\,
\HH^0({{X}}, \Lscr_n)\otimes_k \HH^0({{X}}, \Lscr_m^{\sigma^n})
\,\,\,{\buildrel\phi\over\too}\,\,\,
\HH^0({{X}}, \Lscr_{n+m}) \,\,=\,\, B_{n+m}.
\end{array}
$$

The ring $B$ has two significant properties. First, the way it has been
constructed  ensures that the natural isomorphism $S_1\cong \HH^0(\mathbb
P^2,\Oscr_{\mathbb P^2}(1)) \cong B_1$  induces a ring homomorphism $\phi: S\to
B$. With a little more work using  the Riemann-Roch theorem  one can
even show that  $B\cong S/gS$ for some $g\in S_3$. Secondly---and
this will be explained in more detail in the next section---$\qgr
B\cong \coh(E)$. The latter fact allows one to obtain a detailed  understanding
of the structure of $B$ and the former allows one to pull  this information
back to $S$.

To summarize, the point modules over the Sklyanin algebra $\Skl_3$  are
determined by an automorphism of an elliptic curve $E$ and the geometry of  $E$
allows one to determine the structure of $\Skl_3$.  As is shown in \cite{ATV1}
this technique works more generally and this leads
to the following theorem.

\begin{theorem}\label{ATV-S} {\rm \cite{ATV1,Steph1,Steph2}}
The AS-regular rings $R$ of dimension $3$ are classified.
They are all noetherian domains
with the Hilbert series of a weighted polynomial ring
$k[x,y,z]$; thus the $(x,y,z)$ can
be given  degrees $(a,b,c)$ other than $(1,1,1)$.

Moreover,  $R$ always maps homomorphically onto a twisted homogeneous
coordinate ring $B=B(X,\Lscr,\sigma)$, for some
scheme $X$. Thus $\coh(X) \simeq \qgr B
\hookrightarrow \qgr R$.  \end{theorem}

In this result, Artin, Tate and Van den Bergh \cite{ATV1} classified the
algebras generated in degree one, while Steph\-enson \cite{Steph1,Steph2} did
the general case.

There are strong arguments (see \cite{Bondal} or
\cite[Section~11]{SV})  for saying that the noncommutative analogues of  the
projective plane are precisely
the categories $\qgr R$, where  $R$ is an AS-regular ring
with the Hilbert series $1/(1-t)^3$ of the
unweighted polynomial ring $k[x,y,z]$. So consider this class, which
 clearly includes the Sklyanin
algebra.  The second paragraph of the theorem can now be refined to say that
{\it either} $X=\mathbb P^2$, in which case  $\qgr R\simeq \coh(\mathbb P^2)$,
{\it or} $X=E$ is a cubic curve in $\mathbb P^2$.  Thus, the theorem can be
interpreted as saying that {\it noncommutative projective planes
are either equal to  $\mathbb P^2$ or contain a commutative curve  $E$.}

\section{Twisted homogeneous coordinate rings} \label{sec1.6}
\setzero\vskip-5mm \hspace{5mm}

The ideas from \cite{ATV1} outlined in the last section have had many
other applications, but before we discuss them we need to analyse
twisted homogeneous coordinate rings  in more detail.
The following exercise may give the reader a feel for the construction.

\begin{exercise} \label{twisted-p} \rm Perhaps the simplest algebra appearing in
the theory of quantum groups is the quantum (affine) plane  $k_q[x,y] =
k\{x,y\}/(xy-qyx)$, for $q\in k^*$. Prove that
$k_q[x,y]\cong B(\pone,\Oscr_{\pone}(1), \sigma)$ where $\sigma$ is defined by
$\sigma(a:b) = (a:qb)$, for $(a:b)\in \pone$.  \end{exercise}

For the rest of the section, fix a $k$-scheme $X$ with an invertible  sheaf
$\Lscr$ and  automorphism $\sigma$. When $\sigma=1$, the homogeneous coordinate
ring $B(X,\Lscr)=B(X,\Lscr,1)$ is a standard construction and one has Serre's
fundamental theorem: If $\Lscr$ is ample then $\coh(X)\simeq \qgr(B)$. As was
hinted in the last section, this does  generalize to the noncommutative case,
provided one changes the definition of ampleness. Define $\Lscr$ to be {\it
$\sigma$-ample } if, for all $\Fscr\in\coh(X)$, one has  $\HH^q(X,\Fscr\otimes
\Lscr_n)=0$ for all $q> 0$ and all $n\gg 0$. The na\"\i ve generalization of
Serre's Theorem then holds.

 \begin{theorem} \label{AV} {\rm (Artin-Van den Bergh \cite{AVdB})}
 Let $X$ be a projective scheme with an automorphism $\sigma$ and let
  $\Lscr$ be a $\sigma$-ample invertible sheaf. Then $B=B(X,\Lscr,\sigma)$ is
a  right noetherian cg ring such that $\qgr(B)\simeq \coh(X)$.
\end{theorem}

This begs the question of precisely which line bundles are
$\sigma$-ample.
 A simple application of the Riemann-Roch Theorem shows that
\begin{equation}\label{rroch}
\text{\it if $X$ is a curve, then
 any ample invertible sheaf is $\sigma$-ample,}
 \end{equation}
and the converse holds for irreducible curves. This explains why
Theorem~\ref{AV} could be applied to the factor of the Sklyanin
algebra in the last section.

For higher
dimensional varieties the situation is more subtle and is described by
the following result, for which we need some notation.
Let $X$ be a projective scheme and write  $A^1_{\Num}(X)$ for the set of
Cartier divisors of $X$ modulo  numerical equivalence.
Let $\sigma$ be an automorphism of $X$ and let $P_\sigma$ denote its
induced action on  $A^1_{\Num}(X)$. Since $A^1_{\Num}(X)$ is a finitely
generated  free abelian group,
$P_\sigma$ may be represented by a matrix and
$P_\sigma$ is called {\it quasi-unipotent}\label{quasiunipot-index}
 if all the eigenvalues of this matrix are roots of unity.

\begin{theorem}\label{sigma-ample} {\rm (Keeler \cite{Ke})}
If $\sigma$ be an automorphism of a  projective scheme $X$
then:
\begin{enumerate}
\item[(1)]  $X$ has a $\sigma$-ample line bundle if and only if
$P_\sigma$ is quasi-unipotent. If  $P_\sigma$ is quasi-unipotent, then
all ample line bundles are $\sigma$-ample.
\item[(2)] In Theorem~\ref{AV}, $B$ is also left noetherian.
\end{enumerate}\end{theorem}

There are two comments that should be made about
Theorem~\ref{sigma-ample}. First, it is standard that  $\GKdim B(X,\Lscr) =
1+\dim X$, whenever $\Lscr$ is  ample. However, it can happen that
$\GKdim B(X,\Lscr,\sigma)>1+\dim X$.
Secondly, one can still construct
$B(X,\Lscr,\sigma)$ when $\Lscr$ is ample but $P_\sigma$ is not
quasi-unipotent, but the resulting algebra is rather unpleasant.
Indeed, possibly  after replacing $\Lscr$ by some
$\Lscr^{\otimes n}$, $B(X,\Lscr,\sigma)$  will be a non-noetherian algebra
of  exponential growth.  See \cite{Ke} for the details.

\section{Noncommutative curves and surfaces}\label{sec2}
\setzero\vskip-5mm \hspace{5mm}

 As we have seen, twisted homogeneous coordinate rings are fundamental to the
 study of noncommutative projective planes. However, a more natural starting
 place would be cg algebras of Gelfand-Kirillov dimension two since,
 as we suggested in the introduction, these should correspond to noncommutative
 curves. Their structure is particularly simple.

 \begin{theorem}\label{thm3.3} {\rm \cite{Staf5}}
 Let $R$ be a cg domain of GK-dimension
$2$ generated in degree one.  Then
 there exists an irreducible curve $Y$ with automorphism $\sigma$ and
ample  invertible sheaf $\Lscr$ such that
  $R$ embeds into the twisted homogeneous coordinate  ring
$B(Y,\Lscr,\sigma)$ with finite index. Equivalently,  $R_n\cong \HH^0(Y,
\Lscr\otimes \Lscr^{\sigma}\otimes\cdots\otimes \Lscr^{\sigma^{n-1}})$ for
$ n\gg 0.$
\end{theorem}

By (\ref{rroch})  we
may  apply Theorem~\ref{AV} to obtain part (1) of the next result.

 \begin{corollary}\label{cor3.4} Let $R$ be as in Theorem~\ref{thm3.3}.
Then:
\begin{enumerate}
\item[(1)] $R$ is a noetherian domain with  $\qgr R\simeq\coh(Y)$.
In particular, $\qgr R\simeq \qgr C$ for the commutative ring
$C=B(Y,\Lscr,\Id)$.

\item[(2)] If $|\sigma|<\infty$ then $R$ is a finite module over its centre. If
$|\sigma|=\infty$, then $R$ is a primitive ring with at most
two height one prime ideals.
\end{enumerate}\end{corollary}

If $R$ is not generated in degree one, then the analogue of
Theorem~\ref{thm3.3}  is more subtle, since more complicated algebras appear.
See \cite{Staf5} for the details. One should really make a further
generalization by allowing $R$ to be prime rather than a domain and to allowing
$k$ to be arbitrary (since this allows one to consider the projective analogues
of classical orders  over Dedekind domains).  Theorem~\ref{thm3.3} and
Corollary~\ref{cor3.4} do generalize appropriately
 but the results are more
technical. The details can be found in \cite{Staf6}.

Although these results  are
satisfying they are really only half of the story. As in the commutative
case one would also like to define noncommutative curves
abstractly and then show that they can indeed be described by
graded rings of the appropriate form. Such a result appears in \cite{RV}
but to state it we need a definition.

Let $\Cscr$ be   an $\Ext$-finite abelian category of finite homological
dimension with derived category of bounded complexes $D^b(\Cscr)$.
  Recall that a cohomological functor $H:D^b(\Cscr)\to\mod(k)$ is of
\emph{finite type} if, for $A\in D^b(\Cscr)$, only a finite number of the
$H(A[n])$  are non-zero. The category  $\Cscr$ is
\emph{saturated}\label{saturated-index} if  every cohomological functor
$H:D^b(\Cscr)\to\mod(k)$ of finite type is of the form $\Hom(A,-)$ (that is,
$H$ is representable). If $X$ is a smooth projective scheme, then
$\coh(X)$ is saturated \cite{Bondal4}, so it is not unreasonable to
use this as part of the definition of a ``noncommutative smooth curve.''

\begin{theorem} {\rm (Reiten-Van den Bergh \cite[Theorem~V.1.2]{RV})}
\label{corsat} Assume that $\Cscr$ is a connected
saturated hereditary noetherian category.
   Then  $\Cscr$ has one of the following forms:
\begin{enumerate}
\item[(1)] $\mod(\Lambda)$ where $\Lambda$ is an indecomposable
 finite dimensional  hereditary algebra.
\item[(2)] $\coh(\Oscr)$ where $\Oscr$ is a sheaf of hereditary
$\Oscr_X$-orders over a smooth connected projective curve $X$.
\end{enumerate}
\end{theorem}

It is easy to show that the abelian categories appearing in parts (1) and
(2) of this theorem  are of the form $\qgr R$ for a graded ring $R$
with $\GKdim R\le 2$, and so  this result can be regarded as a
partial converse to Theorem~\ref{thm3.3}. A discussion of the saturation
condition for noncommutative algebras may be found in \cite{BVdB}.

If one accepts that noncommutative projective curves and planes have been
classified, as we have argued, then the natural next step is to attempt to
classify all noncommutative surfaces and this has been a major focus of recent
research.
This program is discussed in detail in \cite[Sections~8--13]{SV} and so here we
will be very brief. For the sake of argument we will assume that an
(irreducible) noncommutative surface is $\qgr R$ for a noetherian cg domain
$R$ with $\GKdim R = 3$, although the precise definition is as yet unclear.
For example, Artin \cite{Ar2} demands that $\qgr R$ should also possess
 a dualizing complex in the sense of Yekutieli \cite{Ye}.
Nevertheless in attempting to classify surfaces it is natural to
mimic the commutative proof:

\begin{itemize}
\item[(a)] Classify noncommutative surfaces up to birational
equivalence; equivalently  classify the associated graded division rings of
fractions for  graded domains $R$ with $\GKdim R=3$.
Artin \cite[Conjecture~4.1]{Ar2} conjectures that
these division rings  are known.

\item[(b)] Prove a version of Zariski's theorem that asserts that one can pass
from any smooth surface to a birationally equivalent one by
successive blowing up and down. Then find minimal models within each
equivalence class.
\end{itemize}

Van den Bergh has created a noncommutative theory of blowing up and down
\cite{VdB19, VdB24} and used this to answer part (b) in a number of special
cases. A key fact in his approach is that (after minor modifications) each
known  example of a noncommutative surface $\qgr R$   contains an embedded
commutative curve $\Cscr$, just as $\qgr(\Skl_3)\hookleftarrow \coh(E)=E$ in
Section~2. This is important since he needs to blow up points on
that subcategory. In general,  define
 {\it a point in $\qgr R$} to be
$\pi(M)$ for a point module  $M\in\gr R$. Given such a point $p$, write
$p=\pi(R/I) =\Rscr/\Iscr$. Mimicking the classical situation we would like
to write
\begin{equation}\label{rees}
\Bscr = \Rscr\oplus \Iscr \oplus \Iscr^2\oplus\cdots,
\end{equation}
and then define the blow-up of $\qgr R$ to be  the
category $\qgr \Bscr$ of  finitely generated
 graded $\Bscr$-modules modulo those
that are right bounded.  However, there are two  problems. A minor one is that
 $\Iscr$ needs to be twisted to take into account the shift functor
on $\qgr R$. The major one is that $I$ is only a one-sided ideal of $R$, and so
there is no natural multiplication on $\Bscr$.
 To circumvent these problems, Van den Bergh \cite{VdB19}
 has to define $\Bscr$ in a more subtle category so that it is indeed
an algebra. It is then
quite hard  to prove that $\qgr \Bscr$   has the appropriate
properties.

\section{Hilbert schemes} \label{sec1.7}
\setzero\vskip-5mm \hspace{5mm}

Since point modules and twisted homogeneous coordinate rings have proved
so useful, it is natural
to ask how  generally these  techniques can be applied. In particular,
one needs to understand when point modules, or other classes of modules, can be
parametrized by a scheme. Indeed, even for point modules over surfaces the
answer was unknown until recently and this is obviously rather important for
the program outlined in the last section.

The best positive result is due to Artin,
Small  and Zhang \cite{ASZ,AZ2}, for which we need a definition. A $k$-algebra
$R$ is called {\it strongly noetherian} if $R\otimes_kC$ is noetherian for all
noetherian commutative $k$-algebras $C$.

\begin{theorem}\label{hilb-schemes}
{\rm (Artin-Zhang  \cite[Theorems~E4.3 and E4.4]{AZ2})}
  Assume that $R$ is a strongly
 noetherian, cg algebra and fix $h(t)=\sum h_it^i\in k[[t]].$
 Let $\Cscr$ denote the set  of cyclic $R$-modules
 $M=R/I$ with  Hilbert series $h_M(t)=\sum \dim_k(M_i)t^i$ equal to $h(t)$.
Then:\begin{enumerate}
\item[(1)]
$\Cscr$ is naturally parametrized by a (commutative) projective scheme.
\item[(2)] There exists an integer $d$ such that, if $M=R/I\in \Cscr$, then
$I$ is generated in degrees $\leq d$ as a right ideal of $R$.
\end{enumerate}\end{theorem}

In particular, if $R$ is  a strongly noetherian cg algebra generated in degree
one, then the set of  point modules is naturally parametrized by a projective
scheme $\Pscr$. In this case one can further show that {\it the shift functor
$M\mapsto M_{\geq 1}[1]$ induces an automorphism $\sigma$ of $\Pscr$.} Thus
one can form the corresponding twisted homogeneous  coordinate rings
$B=B(\Pscr,\Lscr,\sigma)$ and
  for an appropriate line bundle $\Lscr$ there will
exist a  homomorphism $\phi:R\to B$.  Determining when $\phi$
is surjective is probably quite subtle. This result cannot be
used to shorten the arguments about the Sklyanin algebra $\Skl_3$ given in
Section~2, since one needs to use
 $B(E,\Lscr,\sigma)$  to prove that $\Skl_3$ is noetherian.

Although we have concentrated on point modules, more general classes
of modules are also important. An example where {\it line modules} (modules
$M$ with the Hilbert series of $k[x,y]$)
 are needed in a classification problem appears in \cite{vancliff}.

How strong is the strongly noetherian hypothesis? Certainly most of the
standard  examples of noetherian cg algebras (including the Sklyanin algebras)
are strongly noetherian (see \cite[Section~4]{ASZ}) and so one might hope that
this is always the case. But in fact, as Rogalski \cite{Rog} has shown, cg
noetherian algebras that are not strongly noetherian exist in profusion.

These examples are constructed as subrings of
$B=B(\mathbb P^n, \Oscr_{\mathbb P^n}(1), \sigma)$
for an appropriate automorphism $\sigma$. Given
$\sigma\in \Aut(\mathbb P^n)$,
 pick $c\in \mathbb P^n$ and set $\Cscr=\{c_i =\sigma^{-i}(c): i\in
\mathbb N\}$. Then $\Cscr$ is called {\it critically dense} if, for any
infinite subset $\Dscr\subseteq \Cscr$, the Zariski closure  of $\Dscr$ equals
$ \mathbb P^n$. This is not a particularly stringent condition, since it holds
for a generic set of $(\sigma,c) \in \Aut(\mathbb P^n)\times \mathbb P^n.$
 Corresponding to $c$  one has the point
module $M=B/VB$ for some codimension one subspace $V=V(c)\subseteq B_1$.
Rogalski's example is then simply $S(\sigma,c) = k\langle V\rangle \subset B$,
and it has remarkable properties:

\begin{theorem}\label{rog} {\rm (Rogalski \cite{Rog})}
 Keep the above notation. Assume that
 $\sigma\in \Aut(\mathbb P^n)$ and
$c\in \mathbb P^n$ for $n\geq 2$ are such that $\Cscr$ is critically dense.
Then:
\begin{enumerate}
\item[(1)] $S=S(\sigma,c)$ is always noetherian but never strongly noetherian.

\item[(2)] The point modules for $S$ are not naturally parametrized by
 a projective  scheme.

\item[(3)] $S$ satisfies the condition $\chi_1$ but not the condition $\chi_2$,
as defined below.
Moreover, $\qgr S$ has finite cohomological dimension.

\item[(4)]  The category $\qgr S$ is not Ext-finite; indeed if $\Sscr=\pi(S)\in
\qgr S$, then
$\HH^1(\Sscr) = \Ext^1_{\qgr S}(\Sscr,\Sscr)$ is infinite dimensional.
\end{enumerate}\end{theorem}

Some comments about the theorem are in order. First, the point modules for
$S=S(\sigma,c)$  are actually parametrized by an ``infinite blowup of $\mathbb
P^n$'' in the sense that they are parametrized by  $\mathbb P^n$ except that
for each $p\in \Cscr$ one has a  whole family $\Pscr_p$ of point
modules parametrized by $\mathbb P^{n-1}$.
In contrast, the points in $\qgr S$ are
actually parametrized by $\mathbb P^n$ since, if
$M,N\in \Pscr_p$, then $\pi(N)\cong \pi(M)$ in $\qgr S$.

The conditions $\chi_i$ in part (3) are defined as follows: A cg ring $R$
satisfies  $\chi_n$ if, for each $0\leq j\leq n$ and each $M\in \gr R$, one has
$\dim_k \Ext^j_R(k,M)<\infty.$ The significance of $\chi_1$ is that, by
\cite[Theorem~4.5]{AZ}, one can reconstruct $S=S(\sigma,c)$ from $\qgr
S$ and so the  peculiar properties of $S$ are reflected in
$\qgr S$. In particular, part (4) implies that $S$ does not
satisfy $\chi_2$. The  significance  of part (4) is that, for all the
algebras $R$ considered  until now, Serre's finiteness theorem
holds in the sense that $\HH^i(\Fscr)$ is finite
dimensional for all $\Fscr\in\qgr R$ and all $i$.

Here is the simplest example of $S(\sigma,c)$. Pick algebraically independent
elements $p,q\in k$  and
define $\sigma\in \Aut(\mathbb P^2)$ by $\sigma(a \collon b \collon c)
 = (pa\collon qb\collon c)$.
If $c=(1\collon 1\collon 1) \in \mathbb P^2$ then $\Cscr$ is critically dense
and an argument like that of Exercise~\ref{twisted-p} shows that
$$B=k\{x,y,z\}/(zx-pxz,\,zy-qyz,\,yx-pq^{-1}xy)\quad {\rm and} \quad
S(\sigma,c)=k\langle
y-x,\, z-x\rangle.$$
This example was first considered by Jordan \cite{Jo} who
was able to parametrize the point modules for  $S(\sigma,c)$ but was unable to
determine if the ring was noetherian.

Rogalski's examples show that, even for surfaces, the picture is much more
complicated than the discussion of the last section would suggest.
Yet even these examples
  appear in a  geometric framework; indeed they can be
constructed as  blow-ups of $\mathbb P^n$ if one uses the na\"\i ve  approach
of (\ref{rees}).

This works as follows. As before, assume that
$(\sigma,c)\in \Aut(\mathbb P^n)\times \mathbb P^n$ for $n\geq 2$ is
such that $\Cscr$ is critically dense. In
$\coh(\mathbb P^n)$ let $\Iscr_c$ denote the ideal sheaf corresponding to  the
point $c$. If $\Lscr$ is a coherent module over  $\Oscr=\Oscr_{\mathbb P^n}$,
we form a bimodule  $\Lscr_\sigma$ such that as a left module,
$\Lscr_\sigma\cong \Lscr$ but the right action is twisted by $\sigma$: if $s\in
\Lscr_\sigma(U)$ and $a\in  \Oscr_{\mathbb P^n}(\sigma U)$, then $sa\in
\Lscr_\sigma(U) $ is defined by the formula $sa=a^\sigma s$. See
\cite[pp.252-3]{AVdB} for a more formal discussion. Now set
$\Jscr=\Iscr_c\otimes_{\Oscr}\Oscr(1)_\sigma\subseteq  \Oscr(1)_\sigma$ and let
$\Bscr = \Bscr(\sigma,c)=\Oscr \oplus \Jscr \oplus \Jscr^2\oplus\cdots$,  where
$\Jscr^n$ is the image of $\Jscr^{\otimes n}$ in $\Oscr(1)_\sigma^{\otimes
n}\cong  \Oscr(n)_{\sigma^n}$. This does not define a sheaf of rings in the
usual sense since we are ``playing a game of musical chairs with the open sets
\cite[p.252]{AVdB}.'' Nevertheless $\Bscr $ does have an natural graded
algebra structure  and so we can form $\qgr \Bscr$ in the usual way. If
$\sigma=1$ then $\qgr \Bscr$ is simply $coh(X)$, where
$X$ is the blow-up of $\mathbb P^n$ at  $c$. In contrast,
Keeler,  Rogalski and the author  have recently proved:

 \begin{theorem} \label{KRS} {\rm \cite{KRS}}
 Pick $(\sigma,c) \in \Aut(\mathbb P^n)\times
\mathbb P^n$ for $n\geq 2$  such that $\Cscr$ is critically dense.
Then $\Bscr = \Bscr(\sigma,c)$ is noetherian. Moreover
$\qgr(\Bscr) \simeq \qgr S(\sigma,c).$
\end{theorem}

Thus, $\qgr S(\sigma,c)$ is nothing more than the (noncommutative) blow-up of
$\mathbb P^n$ at a point! The differences between this blow-up and Van den
Bergh's are illustrative. Van den Bergh had to work hard to ensure that  the
analogue of the exceptional divisor
really looks like a curve.
 Indeed  much of his formalism is required for just this
reason. In contrast, in Theorem~\ref{KRS} the analogue of the exceptional
divisor (which in this case equals $\Bscr/(\Iscr_{c_{-1}})\Bscr$)
 is actually a
point. This neatly explains the structure  of the  points
in $\qgr S(\sigma,c)$; they are indeed parametrized by $\mathbb P^n$
although the point  corresponding to $c$ (and hence the shifts of this
point, which are nothing more than the  points corresponding to the
$c_i$) are distinguished.

\label{lastpage}


\begin{thebibliography}{aa}

\bibitem{Ar2}
M.~Artin,  Some problems on three-dimensional graded domains,
{\it  Representation theory and algebraic geometry,}
London Math. Soc. Lecture Note
  Ser., vol. 238, Cambridge Univ. Press, Cambridge, 1995, 1--19.


\bibitem{AS}
M.~Artin and W.~Schelter,  Graded algebras of global dimension 3, {\it Adv. in
  Math.},  66 (1987), 171--216.

\bibitem{ASZ} M. Artin, L. W.  Small and J. J.  Zhang,
Generic flatness for strongly noetherian algebras.
{\it J. Algebra}, 221 (1999), 579--610.

\bibitem{Staf5}
M.~Artin and J.~T. Stafford,  Noncommutative graded
domains with quadratic growth, {\it Invent.
  Math.}, 122 (1995), 231--276.

\bibitem{Staf6}
M.~Artin and J.~T. Stafford,   Semiprime graded algebras of
dimension two, {\it  J.~Algebra},  277 (2000), 68--123.


\bibitem{ATV1}
M.~Artin, J.~Tate, and M.~Van~den Bergh,  Some algebras associated to
  automorphisms of elliptic curves, {\it The Grothendieck Festschrift}, vol.~1,
  Birkh\"auser, Boston, 1990, 33--85.

\bibitem{AVdB}
M.~Artin and M.~Van~den Bergh,  Twisted homogeneous coordinate rings, {\it J.
  Algebra}, 133 (1990), 249--271.

 \bibitem{AZ}
M.~Artin and J.~J. Zhang,  Noncommutative projective schemes, {\it Adv. in
 Math.}, 109 (1994), 228--287.


\bibitem{AZ2} M.~Artin and J.~J. Zhang,
  Abstract Hilbert schemes, {\it Algebr. Represent. Theory},
 4 (2001), 305--394.


\bibitem{Bondal4}
A.~I. Bondal and M.~M. Kapranov,  Representable functors, Serre functors,
  and reconstructions,
 {\it  Math. USSR-Izv.}, 35 (1990), 519--541.

\bibitem{Bondal}
A.~I. Bondal and A.~E. Polishchuk, Homological properties of associative
  algebras: the method of helices, {\it Russian Acad. Sci. Izv. Math.}, {\bf 42}
  (1994), 219--260.


\bibitem{BVdB}   A. I. Bondal and M. Van den Bergh, Generators
   and representability of functors in commutative and noncommutative geometry;
   math.AG/0204218 (to appear).


\bibitem{connes}   A. Connes and M. Dubois-Violette,
 Noncommutative finite-dimensional manifolds.
 I. Spherical manifolds and related examples;  math.QA/0107070 (to appear).



\bibitem{Jo} D.~A.~Jordan, The graded algebra generated by two Eulerian
derivatives, {\it  Algebr. Represent. Theory},  4 (2001),
249--275.

\bibitem{Ke}
D.~S. Keeler, Criteria for $\sigma$-ampleness, {\it J. Amer. Math. Soc.},
 13 (2000), 517--532.

\bibitem{KRS} D. S. Keeler, D. Rogalski and J. T. Stafford, work in progress.


\bibitem{KL}
G. R. Krause and T.~H. Lenagan, {\it Growth of algebras and
  Gelfand-Kirillov dimension}, Research Notes in Mathematics, vol. 116,
  Pitman, Boston, 1985.


\bibitem{OF2}
A.~V. Odesskii and B.~L. Feigin,
Sklyanin's elliptic algebras, {\it Functional Anal. Appl.},
23  (1989), no.~3, 207--214.

\bibitem{RV}
I.~Reiten and M.~Van~den Bergh, Noetherian hereditary categories
  satisfying Serre duality,
 {\it J. Amer. Math. Soc.}, 15 (2002), 295--366.


\bibitem{Rog} D.~Rogalski, Examples of generic noncommutative surfaces;
  math.RA/0203180 (to
appear).

\bibitem{rosenberg1}
A.~L. Rosenberg,
 The spectrum of abelian categories and reconstruction of
  schemes,  {\it Rings, Hopf algebras, and Brauer
groups}, Lecture Notes in Pure and Appl. Math., vol. 197, Marcel Dekker,
New York,  1998,  257--274.

  \bibitem{vancliff} B. Shelton and M. Vancliff, Schemes of line modules I,
 {\it  J. London Math. Soc.};
 www.uta.edu/math/vancliff/R/ (to appear).

\bibitem{Skly1}
E.~K. Sklyanin,  Some algebraic structures connected to the
Yang-Baxter
  equation, {\it Functional Anal. Appl.}, 16 (1982), 27--34.


\bibitem{Staf4}
J.~T. Stafford and J.~J. Zhang,  Examples in noncommutative projective
  geometry, {\it Math. Proc. Cambridge Philos. Soc.}, 116 (1994), 415--433.


\bibitem{SV} J.~T. Stafford and  M.~Van~den Bergh, Noncommutative curves
and noncommutative surfaces, {\it Bull. Amer. Math. Soc.},   38
(2001), 171--216.

\bibitem{Steph1} D.~R. Stephenson, Artin-Schelter regular algebras of
global dimension three, {\it J. Algebra}, 183 (1996), 55--73.

\bibitem{Steph2} D.~R. Stephenson, Algebras associated to elliptic curves, {\it
Trans. Amer. Math. Soc.},  349 (1997), 2317--2340.


\bibitem{VdB19}
M.~Van~den Bergh, Blowing up of noncommutative smooth surfaces,
{\it Mem. Amer. Math. Soc.}, 154 (2001), no. 734.

 \bibitem{VdB24}
 M.~Van~den Bergh,  Abstract blowing down,
  {\it Proc. Amer. Math. Soc.}, 128 (2000),  375--381.


\bibitem{Ye}
A.~Yekutieli,  Dualizing complexes over noncommutative graded
algebras, {\it J.  Algebra},  153 (1992), 41--84.

\end{thebibliography}
\end{document}